# Escape of mass in zero-range processes with random rates

**Pablo A. Ferrari [1],[*] and Valentin V. Sisko [2],[†]**

*Universidade de São Paulo*

**Abstract:** We consider zero-range processes in $\mathbb{Z}^d$ with site dependent jump rates. The rate for a particle jump from site $x$ to $y$ in $\mathbb{Z}^d$ is given by $\lambda_x g(k) p(y-x)$, where $p(\cdot)$ is a probability in $\mathbb{Z}^d$, $g(k)$ is a bounded nondecreasing function of the number $k$ of particles in $x$ and $\lambda = \{\lambda_x\}$ is a collection of i.i.d. random variables with values in $(c, 1]$, for some $c > 0$. For almost every realization of the environment $\lambda$ the zero-range process has product invariant measures $\{\nu_{\lambda,v} : 0 \leq v \leq c\}$ parametrized by $v$, the average total jump rate from any given site. The density of a measure, defined by the asymptotic average number of particles per site, is an increasing function of $v$. There exists a product invariant measure $\nu_{\lambda,c}$, with maximal density. Let $\mu$ be a probability measure concentrating mass on configurations whose number of particles at site $x$ grows less than exponentially with $\|x\|$. Denoting by $S_\lambda(t)$ the semigroup of the process, we prove that all weak limits of $\{\mu S_\lambda(t), t \geq 0\}$ as $t \to \infty$ are dominated, in the natural partial order, by $\nu_{\lambda,c}$. In particular, if $\mu$ dominates $\nu_{\lambda,c}$, then $\mu S_\lambda(t)$ converges to $\nu_{\lambda,c}$. The result is particularly striking when the maximal density is finite and the initial measure has a density above the maximal.

## 1. Introduction

In the zero-range process there are a finite number of particles at each site of $\mathbb{Z}^d$. At a rate depending monotonically on the number of particles at the site, one of the particles jumps to another site chosen independently with a transition probability function. The rate at which particles leave any site is bounded. When the rate at each site $x$ is multiplied by a random variable $\lambda_x$ chosen at time zero independently of the process, the system may show a phase transition in the density. For almost every realization of the environment $\lambda$ the zero-range process has product invariant measures $\{\nu_{\lambda,v} : 0 \leq v \leq c\}$ parametrized by $v$, the average total jump rate from any given site. The density of a measure is the asymptotic number of particles per site (when this exists). For each $v \leq c$ the invariant measure $\nu_{\lambda,v}$ has density $\rho(v)$, which is an increasing function of $v$. Our main result is to start the system with a measure concentrating mass in configurations not growing too fast (see (3) below) and show that the distribution of the process as time goes to infinity is dominated by the maximal measure $\nu_{\lambda,c}$. This is particularly interesting when $\rho(c) < \infty$ and the initial density of $\mu$ is strictly bigger than $\rho(c)$. In this case we say that there is an "escape of mass". When the initial distribution dominates the maximal invariant measure, the process converges to the maximal invariant measure.

---

[*]Supported in part by FAPESP.
[†]Supported by FAPESP (2003/00847–1) and CNPq (152510/2006–0).
[1]Departamento de Estatística, Instituto de Matemática e Estatística, Universidade de São Paulo, Caixa Postal 66281, CEP 05311–970 São Paulo, SP, Brazil, e-mail: pablo@ime.usp.br, url: www.ime.usp.br/~pablo
[2]IMPA, Estrada Dona Castorina 110, CEP 22460-320 Rio de Janeiro, Brasil, e-mail: valentin@impa.br, url: www.ime.usp.br/~valentin

The zero-range process appeared first as a network of queues when Jackson [13] showed that the product measures are invariant for the process in a finite number of sites. Spitzer [22] introduced the process in a countable number of sites as a model of infinite particle system with interactions. The existence of the process has been proved by Holley [12] and Liggett [17, 19]. We use Harris [11] direct probabilistic construction which permits the particles to be distinguishable, so one can follow the behavior of any particular particle. Using Liggett's [18] approach, Andjel [1] gave a description of the set of invariant measures for the zero-range process in some cases. Balázs, Rassoul-Agha, Seppäläinen, and Sethuraman [4] studied the case of rates bounded by an exponential function of $k$ in a one dimensional asymmetric model.

The study of conservative interacting particle systems in random environment was proposed simultaneously by Benjamini, Ferrari and Landim [5] and Evans [7], who observed the existence of phase transition in these models; see also Krug and Ferrari [15]. Benjamini, Ferrari and Landim [5], Krug and Seppäläinen [20] and Koukkous [14] investigated the hydrodynamic behavior of conservative processes in random environments; Landim [16] and Bahadoran [3] considered the same problem for non-homogeneous asymmetric attractive processes; Gielis, Koukkous and Landim [9] deduced the equilibrium fluctuations of a symmetric zero-range process in a random environment; Andjel, Ferrari, Guiol and Landim [2] proved the convergence to the maximal invariant measure for a one-dimensional totally asymmetric nearest-neighbor zero-range process with random rates. This phenomenon is studied by Seppäläinen, Grigorescu and Kang [10] in one dimension. Evans and Hanney [8] have recently published a review paper on the zero-range process which includes many references to the mathematical physics literature.

Section 2 includes definitions, results and at the end a summary of the contents of the other sections.

## 2. Results

We study the zero-range process with site dependent jump rates. Let $\mathbb{N} = \{0, 1, 2, \dots\}$ and give $\mathbb{N}$ the discrete topology. It would seem natural to take $\mathcal{X} = \mathbb{N}^{\mathbb{Z}^d}$ for the state space, but for topological reasons, let us begin by setting

$$\overline{\mathbb{N}} = \mathbb{N} \cup \{\infty\}.$$

We give $\overline{\mathbb{N}}$ the topology of one point compactification and take $\overline{\mathcal{X}} = \overline{\mathbb{N}}^{\mathbb{Z}^d}$ with the product topology for the state space. The set $\overline{\mathcal{X}}$ is compact. We associate with $\overline{\mathcal{X}}$ the Borel $\sigma$-field. The product topology on $\overline{\mathcal{X}}$ is metrizable. For $x = (x_1, \dots, x_d) \in \mathbb{Z}^d$, denote the sup-norm of $x$ by

$$\|x\| = \max_{i=1,\dots,d} |x_i|.$$

Let $\gamma : \overline{\mathbb{N}} \to [0, 2]$ be such that $\gamma(0) = 2$, $\gamma(n) = 1/n$, $n = 1, 2, \dots$, and $\gamma(\infty) = 0$. For instance, the metric

$$d(\eta, \xi) = \sum_{x \in \mathbb{Z}^d} \frac{1}{2^{\|x\|}} \big|\gamma(\eta(x)) - \gamma(\xi(x))\big|$$

is compatible with the product topology on $\overline{\mathcal{X}}$. The set $\overline{\mathcal{X}}$ is a complete separable metric space.

Fix $0 < c < 1$ and consider a collection $\lambda = \{\lambda_x\}_{x \in \mathbb{Z}^d}$ taking values in $(c, 1]$ such that $c = \inf_{x \in \mathbb{Z}^d} \lambda_x$. We call $\lambda$ the *environment*. Let $p : \mathbb{Z}^d \to [0, 1]$ be a probability on $\mathbb{Z}^d$: $\sum_{x \in \mathbb{Z}^d} p(x) = 1$. We assume that the range of $p$ is bounded by some $M > 0$: $p(x) = 0$ if $\|x\| > M$. Moreover, suppose that the random walk with transition function $p(x, y) = p(y - x)$ is irreducible.

Let $g : \overline{\mathbb{N}} \to [0, 1]$ be a nondecreasing continuous function with $0 = g(0) < g(1)$ and $g(\infty) = \lim g(k) = 1$.

The zero-range process in the environment $\lambda$ is a Markov process informally described as follows. Initially distribute particles on the lattice $\mathbb{Z}^d$, then if there are $k$ particles at site $x$, at rate $\lambda_x g(k) p(y - x)$ a particle jumps from $x$ to $y$. In Section 5 we recall the construction of a process $\eta_t$ with this behavior as a function of a Poisson process in $\mathbb{Z}^d \times \mathbb{R}$, *à la Harris*. Let $\{S_\lambda(t), t \geq 0\}$ be the semigroup associated to this process, that is,

$$S_\lambda(t) f(\eta) = \mathbb{E}[f(\eta_t) \mid \eta_0 = \eta].$$

where $\mathbb{E}$ is expectation and $\eta_t = \eta_{\lambda_t}$ is the process with fixed environment $\lambda$. The corresponding generator $L_\lambda$, defined by

$$L_\lambda f(\eta) = \frac{d}{dt} S_\lambda(t) f(\eta) \Big|_{t=0},$$

acts on cylinder continuous functions $f : \overline{\mathbb{N}}^{\mathbb{Z}^d} \to \mathbb{R}$ as follows:

$$(L_\lambda f)(\eta) = \sum_{x \in \mathbb{Z}^d} \sum_{y \in \mathbb{Z}^d} \lambda_x \, p(y - x) \, g(\eta(x)) \, [f(\eta^{x,y}) - f(\eta)].$$

where $\eta^{x,y} = \eta - \delta_x + \delta_y$ and $\delta_z \in \mathcal{X}$ is the configuration with just one particle at $z$ and no particles elsewhere; addition of configurations is performed componentwise. We set $\infty \pm 1 = \infty$.

The natural state space for this Markov process is $\mathcal{X}$ rather than $\overline{\mathcal{X}}$. From the construction *à la Harris* it is possible to see that if the standard Markov process whose semigroup is given by $S_\lambda(t)$ is started in $\mathcal{X}$, then it never leaves $\mathcal{X}$: if $\mu(\mathcal{X}) = 1$, then $\mu S_\lambda(t)(\mathcal{X}) = 1$ for any $t$.

For each $v \in [0, c]$ and environment $\lambda$, denote $\nu_{\lambda,v}$ the product measure with marginals

(1) $$\nu_{\lambda,v}\{\xi : \xi(x) = k\} = \frac{1}{Z(v/\lambda_x)} \frac{(v/\lambda_x)^k}{g(k)!},$$

where we use the notation $g(k)! = g(1) \cdots g(k)$ and $g(0)! = 1$;

(2) $$Z(u) = \sum_{k \geq 0} \frac{u^k}{g(k)!}$$

is the normalizing constant. These measures are invariant for the process [1, 13, 22]. In some cases it is known that all invariant measures (concentrated on $\mathcal{X}$) are convex combinations of measures in $\{\nu_{\lambda,v} : 0 \leq v \leq c\}$ (see [1, 2]).

To define the standard partial order for probability measures on $\overline{\mathcal{X}}$ let $\eta \leq \xi$ if $\eta(x) \leq \xi(x)$ for all $x \in \mathbb{Z}^d$. A real valued function $f$ defined on $\overline{\mathcal{X}}$ is *increasing* if $\eta \leq \xi$ implies that $f(\eta) \leq f(\xi)$. If $\mu$ and $\nu$ are two probability measures on $\overline{\mathcal{X}}$, $\mu \leq \nu$ if $\int f d\mu \leq \int f d\nu$ for all increasing continuous functions $f$. In this case we say that $\nu$ *dominates* $\mu$. This is equivalent to the existence of a probability measure $\bar{\nu}$ on $\overline{\mathcal{X}} \times \overline{\mathcal{X}}$ with marginals $\mu$ and $\nu$ such that

$$\bar{\nu}\{(\eta, \xi) : \eta \leq \xi\} = 1,$$

(coupling); see Theorem 2.4 of Chapter II in [19].

Since $\overline{\mathcal{X}}$ is compact, any sequence of probability measures on $\overline{\mathcal{X}}$ is tight, and therefore, has a weakly convergent subsequence.

Our main theorem holds for measures $\mu$ on $\overline{\mathcal{X}}$ giving total mass to configurations for which the number of particles in $x$ increases less than exponentially with $\|x\|$. That is, measures satisfying

(3) $$\sum_{n=1}^{\infty} e^{-\beta n} \sum_{x : \|x\| = n} \eta(x) < \infty \quad \mu\text{-a.s. for all } \beta > 0.$$

The product measure $\nu_{\lambda, v}$ obviously satisfies (3).

We consider random rates $\lambda = \{\lambda_x\}_{x \in \mathbb{Z}^d}$, a collection of independent identically distributed random variables in $(c, 1]$. Call $\mathcal{P}$ and $\mathcal{E}$ the probability and expectation induced by these variables. Assume that for any $\varepsilon > 0$, $\mathcal{P}(\lambda_0 \in (c, c + \varepsilon)) > 0$.

**Theorem 1.** *Let $\mu$ be a probability measure on $\overline{\mathcal{X}}$ satisfying (3). Then $\mathcal{P}$-a.s.*

(i) *Every weak limit of $\mu S_\lambda(t)$ as $t$ tends to infinity is dominated by $\nu_{\lambda, c}$.*
(ii) *If $\nu_{\lambda, c} \leq \mu$ then $\mu S_\lambda(t)$ converges to $\nu_{\lambda, c}$ as $t$ goes to infinity.*

The result is better understood using the notion of density of particles. Recall that $\lim g(k) = 1$ and notice that the function $Z : [0, 1) \to [0, \infty)$ defined in (2) is analytic. Let $R : [0, 1) \to [0, \infty)$ be the strictly increasing function defined by

$$R(u) = \frac{1}{Z(u)} \sum_{k \geq 0} k \frac{u^k}{g(k)!} = u \frac{Z'(u)}{Z(u)}.$$

It is easy to see that $R$ is onto $[0, \infty)$. Under the measure $\nu_{\lambda, v}$ the expected number of particles (density) at site $x$ is

(4) $$\nu_{\lambda, v}[\eta(x)] = R(v/\lambda_x),$$

and the expected value of the jump rate is

$$\nu_{\lambda, v}[\lambda_x g(\eta(x))] = v.$$

Since $v/\lambda_x < 1$, for any $v \in [0, c]$ and $x$,

(5) $$\nu_{\lambda, c}[\eta(x)] = \lim_{v \to c} R(v/\lambda_x) < \infty.$$

Since the rate distribution is translation invariant, taking the average with respect to the rates, the mean number of particles per site is

$$\rho(v) := \int \mathcal{P}(d\lambda_0) R(v/\lambda_0).$$

For $v \in [0, c)$, $\rho(v) < \infty$. Depending on the distribution of $\lambda_0$, two cases are possible: $\rho(c) < \infty$ and $\rho(c) = \infty$. Since $R(u)$ is a nondecreasing nonnegative function,

(6) $$\rho(c) = \lim_{v \nearrow c} \rho(v),$$

The equation also holds when $\rho(c) = \infty$.

For $v \in [0, c]$, denote $m_v := \int \mathcal{P}(d\lambda) \, \nu_{\lambda,v}$ the measure obtained by first $\mathcal{P}$-choosing an environment $\lambda$ and then choosing a configuration $\eta$ with $\nu_{\lambda,v}$. Under this law $\{\eta(x)\}_{x \in \mathbb{Z}^d}$ are independent identically distributed random variables with average number of particles per site given by $m_v[\eta(0)] = \rho(v)$. By the strong law of large numbers,

$$\text{(7)} \qquad \lim_{n \to \infty} \frac{1}{(2n+1)^d} \sum_{\|x\| \leq n} \eta(x) = \rho(v) \quad m_v\text{-a.s.}$$

Thus, $\mathcal{P}$-a.s., the limit (7) holds $\nu_{\lambda,v}$-a.s.; it also holds when $\rho(v) = \infty$.

For $\eta \in \mathcal{X}$, the *lower asymptotic density* of $\eta$ is defined by

$$\text{(8)} \qquad \underline{D}(\eta) := \liminf_{n \to \infty} \frac{1}{(2n+1)^d} \sum_{\|x\| \leq n} \eta(x),$$

and the *upper asymptotic density* of $\eta$ is defined by

$$\text{(9)} \qquad \overline{D}(\eta) := \limsup_{n \to \infty} \frac{1}{(2n+1)^d} \sum_{\|x\| \leq n} \eta(x).$$

Take some probability measure $\mu$ satisfying (3) and some environment $\lambda$. Let $\tilde{\mu}$ be a weak limit of $\mu S_\lambda(t)$ along a convergent subsequence. Then Theorem 1 (i) implies

$$\text{(10)} \qquad \overline{D}(\eta) \leq \rho(c) \quad \tilde{\mu}\text{-a.s.}$$

Suppose that $\rho(c) < \infty$ and $\mu$ concentrates mass on configurations with lower asymptotic density strictly bigger than $\rho(c)$, that is,

$$\text{(11)} \qquad \underline{D}(\eta) > \rho(c) \quad \mu\text{-a.s.}$$

Inequality (10) says that weak limits of $\mu S_\lambda(t)$ are concentrated on configurations with the upper asymptotic density of $\eta$ not greater than $\rho(c)$. This behavior is remarkable as the process is *conservative*, i.e., the total number of particles is conserved, but in the above limit there is an "escape of mass". Heuristically, a fraction of the particles get stacked at further and further sites with lower and lower rates.

**Sketch of proof.** The proof is based on the study of a family of zero-range processes indexed with $\alpha > 0$; we call them the $\alpha$-truncated process. The $\alpha$-truncated process behaves as the original process but at all times there are infinitely many particles in sites $x$ with $\lambda(x) \leq c + \alpha$. The measure $\nu_\lambda^\alpha$ is invariant for the process. Let the measure $\mu^\alpha$ be the law of a configuration chosen with $\mu$ modified by putting infinitely many particles in sites $x$ with $\lambda(x) \leq c + \alpha$ and leaving the other sites unchanged. We use the fact that there is a density of sites with infinitely many particles to show that the $\alpha$-truncated process starting with $\mu^\alpha$ for $\mu$ satisfying (3) converges weakly to $\nu_\lambda^\alpha$. We prove the convergence using coupling arguments. Two $\alpha$-truncated processes starting respectively with $\mu^\alpha$ and the invariant law $\nu_\lambda^\alpha$ are jointly realized using the so called "basic coupling" [19] which amounts to use the same Poisson processes to construct both marginals. The coupling induces first and second class particles, the last represent the discrepancies between both marginals.

A key element of the proof is the study of the motion of a single tagged second class particle in the $\alpha$-truncated process. The skeleton of the trajectory of each particle is a simple random walk with jump probabilities $p(\cdot)$ absorbed at sites $x$ with $\lambda(x) \leq c + \alpha$. The interaction with the other particles and with the environment $\lambda$ governs the waiting times between jumps but does not affect the skeleton of the motion. We show that with probability one (a) only a finite number of second class particles will visit any fixed site $x$: particles starting sufficiently far away will be absorbed before arriving to $x$ and (b) the finite number of particles hitting $x$ will be eventually absorbed. The weak convergence and the uniqueness of the invariant measure for the $\alpha$-process is a consequence of this result. The $\alpha$-process dominates stochastically the original process (which corresponds to $\alpha = 0$) when both start with the same configuration. Since $\nu_\lambda^\alpha$ converges to the maximal invariant measure as $\alpha \to 0$, this will conclude the proof.

In Section 3 we introduce the $\alpha$-truncated process, and state the two main results which lead to the proof of Theorem 1: the ergodicity of the $\alpha$-truncated process and the fact that it dominates the original process. In the same section we prove Theorem 1. In Section 4 we prove results for the random walk absorbed at sites $x$ with $\lambda_x \leq c + \alpha$, and in Section 5 we graphically construct the process, introduce the relevant couplings and prove the ergodicity and domination results.

## 3. The $\alpha$-truncated process

We introduce a family of zero-range process with infinite number of particles at sites with sufficiently slow rates. Let $\alpha > 0$, $c^\alpha = c + \alpha$ and $\lambda^\alpha = \{\lambda_x^\alpha\}_{x \in \mathbb{Z}^d}$ the truncation given by

$$\lambda_x^\alpha = \begin{cases} c^\alpha & \text{if } \lambda_x \leq c^\alpha, \\ \lambda_x & \text{if } \lambda_x > c^\alpha. \end{cases}$$

For each $\alpha \geq 0$ consider a $\overline{\mathcal{X}}$-valued zero-range process $\eta_t^\alpha$ in the environment $\lambda^\alpha$. We call it the $\alpha$-truncated process or just the truncated process when $\alpha$ is clear from the context. When $\alpha = 0$ we have the original process: $\eta_t^0 = \eta_t$. Partition $\mathbb{Z}^d = \Lambda(\lambda, \alpha) \cup \Lambda^c(\lambda, \alpha)$ with

$$\Lambda(\lambda, \alpha) = \{x \in \mathbb{Z}^d : \lambda_x > c + \alpha\} \quad \text{and} \quad \Lambda^c(\lambda, \alpha) = \{x \in \mathbb{Z}^d : \lambda_x \leq c + \alpha\}.$$

We impose that $\eta_t^\alpha(x) = \infty$ for all $t$ for all $x \in \Lambda^c(\lambda, \alpha)$. The truncated process $\eta_t^\alpha$ is defined in the same way as $\eta_t$ from Section 2 with the following differences. Particles jump as before to $\Lambda^c(\lambda, \alpha)$, but since there are infinitely many particles in $\Lambda^c(\lambda, \alpha)$, the rate of jump from $x \in \Lambda^c(\lambda, \alpha)$ to $y$ is $(c+\alpha)g(\infty)p(y-x)$. Since the number of particles in $x$ is always infinity, this jumps can be interpreted as creation of particles in $y$. Hence the process $\eta_t^\alpha$ can be thought of as evolving in $\mathcal{X}^\alpha := \mathbb{N}^{\Lambda(\lambda, \alpha)}$ with boundary conditions "infinitely many particles at sites in $\Lambda^c(\lambda, \alpha)$".

Let $L_\lambda^\alpha$ be the generator of the $\alpha$-truncated process $\eta_t^\alpha$ and $\{S_\lambda^\alpha(t), t \geq 0\}$ be the semigroup associated to the generator $L_\lambda^\alpha$. We construct this process *à la Harris* in Section 5.

We consider measures on configurations of the processes $\eta_t$ and $\eta_t^\alpha$ as measures on $\overline{\mathcal{X}}$. The product measure $\nu_\lambda^\alpha$ with marginals

$$\nu_\lambda^\alpha\{\xi : \xi(x) = k\} = \begin{cases} \dfrac{1}{Z(c^\alpha/\lambda_x^\alpha)} \dfrac{(c^\alpha/\lambda_x^\alpha)^k}{g(k)!} & \text{if } x \in \Lambda(\lambda, \alpha), \\ \mathbf{1}\{k = \infty\} & \text{if } x \in \Lambda^c(\lambda, \alpha), \end{cases}$$

is invariant for the process $\eta_t^\alpha$. Since $c^\alpha \to c$ and $\lambda^\alpha(x) \to \lambda(x)$ as $\alpha$ goes to zero,

$$(12) \qquad \lim_{\alpha \to 0} \nu_\lambda^\alpha = \nu_{\lambda,c} \quad \text{weakly.}$$

Let $T^\alpha : \overline{\mathcal{X}} \to \overline{\mathcal{X}}$ be the truncation operator defined by

$$(13) \qquad T^\alpha \eta(x) = \begin{cases} \eta(x) & \text{if } \lambda_x > c + \alpha, \\ \infty & \text{if } \lambda_x \leq c + \alpha. \end{cases}$$

The operator $T^\alpha$ induces an operator on measures that we also call $T^\alpha$. Define $\mu^\alpha := T^\alpha \mu$. We clearly have

$$(14) \qquad \mu \leq \mu^\alpha.$$

This domination is preserved by the respective processes:

**Lemma 1.** *Let $\alpha > 0$ and $t \geq 0$. Then $\mu S_\lambda(t) \leq \mu^\alpha S_\lambda^\alpha(t)$.*

The truncated process converges to the invariant measure:

**Proposition 1.** *Let $\mu$ be a probability measure on $\overline{\mathcal{X}}$ satisfying (3). Then for any $\alpha > 0$,*

$$(15) \qquad \lim_{t \to \infty} \mu^\alpha S_\lambda^\alpha(t) = \nu_\lambda^\alpha \quad \mathcal{P}\text{-a.s.}$$

We prove Lemma 1 and Proposition 1 in Section 5.

*Proof of Theorem 1.* For any $\alpha > 0$, Lemma 1 and Proposition 1 imply

$$\limsup_{t \to \infty} \mu S_\lambda(t) \leq \limsup_{t \to \infty} \mu^\alpha S_\lambda^\alpha(t) = \nu_\lambda^\alpha.$$

Item (i) follows by taking $\alpha \to 0$ and applying (12).

To prove item (ii), take $\mu$ such that $\nu_{\lambda,c} \leq \mu$. In the same way as in the proof of Lemma 1, it is easy to see that the semigroup $S_\lambda(t)$, acting on measures, preserves the ordering: $\nu_{\lambda,c} S_\lambda(t) \leq \mu S_\lambda(t)$ for any $t$. Since $\nu_{\lambda,c}$ is invariant, $\nu_{\lambda,c} = \nu_{\lambda,c} S_\lambda(t)$. Therefore, by item (i),

$$\nu_{\lambda,c} = \limsup_{t \to \infty} \nu_{\lambda,c} S_\lambda(t) \leq \limsup_{t \to \infty} \mu S_\lambda(t) \leq \nu_{\lambda,c}. \qquad \square$$

Our task is to prove Proposition 1. The point is that the skeleton of each particle is just a discrete-time random walk with absorption at the sites where $\lambda_x \leq c + \alpha$. Since there is a positive density of those sites, only a finite number of particles will arrive at any fixed finite region. On the other hand, the absorbing sites create new particles. We couple the process with initial measure $\mu^\alpha$ with the process with initial invariant measure $\nu_\lambda^\alpha$ in such a way that new particles are created at the same time in the same sites to both processes. New created particles jump together at both marginals. We show that as time goes to infinity, in both processes only new particles will be present in any finite region.

## 4. Family of independent random walks

Fix $\eta$ such that the inequality in (3) holds. Fix $\alpha > 0$. Since $\eta$ and $\alpha$ are fixed, we omit them in the notation when it is possible. For example, $\Lambda^c(\lambda) := \Lambda^c(\lambda, \alpha)$.

For each $x \in \mathbb{Z}^d$, enumerate the $\eta(x)$ particles at site $x$ in some way and let

$$\zeta = \{\zeta_n(x, i) : x \in \mathbb{Z}^d, i \in \mathbb{N} \cap [1, \eta(x)]\}$$

be a family of independent *discrete-time* random walks with starting points $\zeta_0(x, i) = x$, $x \in \mathbb{Z}^d$, $i \in \mathbb{N} \cap [1, \eta(x)]$ and transitions governed by $p(\cdot)$. We use the notation $\mathbf{P}$ and $\mathbf{E}$ for the law and expectation induced by $\zeta$. Recall $\mathcal{P}$ and $\mathcal{E}$ are the law and expectation induced by the environment $\lambda$. By $\mathcal{P} \times \mathbf{P}$ denote the product measure with marginals $\mathcal{P}$ and $\mathbf{P}$.

For each $(x, i)$ and for each subset $A$ of $\mathbb{Z}^d$, denote

$$\tau(x, i; A) = \min\{n \geq 0 : \zeta_n(x, i) \in A\}$$

the first time the walk hits the set $A$ (this could be $\infty$).

Let us prove that if we consider the random walks in time $[0, \tau(x, i; \Lambda^c(\lambda))]$ only a finite number of walks visit the origin and the number of visits of the origin by each of the walks is finite. More formally, by $N(\lambda, \zeta)$ denote the last time any walk visits the origin before entering in $\Lambda^c(\lambda)$:

$$N(\lambda, \zeta) = \sup \bigcup_x \bigcup_i \{m : m \in [0, \tau(x, y; \Lambda^c(\lambda))] \text{ and } \zeta_m(x, i) = 0\}.$$

**Proposition 2.**

(16) $$(\mathcal{P} \times \mathbf{P})\{(\lambda, \zeta) : N(\lambda, \zeta) < \infty\} = 1.$$

*Proof.* Denote $\theta = \mathcal{P}(\lambda_0 \leq c + \alpha)$. If $\alpha$ is small enough, then $0 < \theta < 1$. Call $E_{x,i}$ the subset of $\mathbb{Z}^d$ visited by the walk $\zeta_n(x, i)$ in the time interval $[0, \tau(x, i; \Lambda^c(\lambda))]$ and denote

$$C_{x,i,N} = \{(\lambda, \zeta) : |E_{x,i}| \geq N\}$$

where $N \geq 0$ and $|E_{x,i}|$ is the number of elements in the set $E_{x,i}$. Since each site of $E_{x,i}$ has probability $\theta$ to be in the set $\Lambda^c(\lambda)$,

(17) $$(\mathcal{P} \times \mathbf{P})(C_{x,i,N}) \leq (1 - \theta)^N \to 0 \quad \text{as } N \to \infty.$$

By hypothesis the random walk with transitions governed by $p(\cdot)$ is irreducible, hence it cannot be confined to a finite region. This implies that the number of new sites visited by time $n$ goes to infinity as $n$ increases. This and (17) implies that

(18) $$(\mathcal{P} \times \mathbf{P})\left(\bigcap_{(x,i):i \leq \eta(x)} \{\tau(x, i; \Lambda^c(\lambda)) < \infty\}\right) = 1.$$

Define

$$D_{x,i} = \{(\lambda, \zeta) : \tau(x, i; \{0\}) < \tau(x, i; \Lambda^c(\lambda))\}.$$

Since the range of the random walk is $M < \infty$, we see that the random walk $\zeta_n(x, i)$ visits at least (the integer part of) $\|x\|/M$ different sites before it reaches the origin. Therefore,

(19) $$(\mathcal{P} \times \mathbf{P})(D_{x,i}) \leq (1 - \theta)^{\|x\|/M}.$$

Thus

$$\sum_{(x,i):i\leq\eta(x)} (\mathcal{P}\times\mathbf{P})(D_{x,i}) \leq \sum_k (1-\theta)^{k/M} \sum_{x:\|x\|=k} \eta(x) < \infty$$

because we assumed $\eta$ satisfies (3). Borel-Cantelli then implies that with $(\mathcal{P}\times\mathbf{P})$ probability one only a finite number of events $D_{x,i}$ happen. Thus, if we consider the random walks in time $[0,\tau(x,i;\Lambda^c(\lambda))]$, then only a finite number of walks visit the origin, and by (18), each walk visits the origin a finite number of times. □

## 5. Construction and coupling

We construct *à la Harris* a Markov process $\eta_t$ on $\overline{\mathcal{X}} = \overline{\mathbb{N}}^{\mathbb{Z}^d}$ corresponding to the above description. Let $(N_{x,y},\ x,y\in\mathbb{Z}^d)$ be a collection of independent Poisson process such that $N_{x,y}$ has intensity $p(y-x)$. If a Poisson event $s$ belongs to a Poisson process $N_{x,y}$, then we say that the event has the origin $x$ and the end $y$. To tune the rate with the environment and the number of particles, we associate to each Poisson event $s\in\cup_{x,y}N_{x,y}$ a random variable $U(s)$, uniform in $[0,1]$, independent of the Poisson processes and independent of the other uniform variables. Since the probability that any two Poisson events from $\cup_{x,y}N_{x,y}$ happen at the same time is zero, all the Poisson events can be indexed by their times, in other words, they can be ordered by their time of occurrence.

The evolution of the process $\eta_t = \eta_{\lambda_t}$ in the environment $\lambda$ is given by the following (deterministic) rule: if the Poisson process $N_{x,y}$ has an event at time $s$ and

(20) $$U(s) < \lambda_x g(\eta_{s-}(x)),$$

then one particle is moved from $x$ to $y$ at that time. Since $g(0)=0$, if no particle is in $x$, then the Poisson event produces no effect in the process in this case.

Using that $p$ is finite range, a percolation argument shows that, for $h$ sufficiently small, $\mathbb{Z}^d$ can be partitioned in finite (random) subsets with the following property: all Poisson events in the interval $[0,h]$ have the origin and the end in the same subset. Since there is a finite number of Poisson events in time interval $[0,h]$ in each of the subsets, the Poisson events can be well ordered by their time of occurrence and the value of $\eta_h$ for each subset can be obtained with the rule (20) proceeding from the first event to the last in each subset. Starting at $\eta_h$, we repeat the construction in the interval $[h,2h]$ and so on. Thus, for any $t$, the process $\eta_t$ is well defined as a function of the Poisson processes and the uniform random variables.

The $\alpha$-truncated process $\eta_t^\alpha$ in the same environment $\lambda$ is also realized as a function of the Poisson processes and uniform variables with a similar rule: if the Poisson process $N_{x,y}$ has an event at time $s$ and

(21) $$U(s) < \lambda_x^\alpha g(\eta_{s-}^\alpha(x)),$$

then one particle is moved from $x$ to $y$ at that time. Rules (20) and (21) induce a natural coupling between the processes $\eta_t$ and $\eta_t^\alpha$. This is the key of the proof of Lemma 1.

We use the notation $\mathbb{P}$ and $\mathbb{E}$ for the probability and expectation induced by the Poisson processes and corresponding uniform associated random variables. Notice that this alea does not depend on $\lambda$.

*Proof of Lemma 1.* Fix a configuration $\eta_0$ and an environment $\lambda$ and let $\eta_0^\alpha(x) = \eta_0(x)$ if $x \in \Lambda(\lambda, \alpha)$ and $\eta_0^\alpha(x) = \infty$ if $x \in \Lambda^c(\lambda, \alpha)$. Let $(\eta_t, \eta_t^\alpha)$ be the coupling obtained by constructing each marginal as a function of the Poisson processes $(N_{x,y}, x, y \in \mathbb{Z}^d)$ and uniform random variables $(U(s), s \in \cup_{x,y} N_{x,y})$ following rules (20) and (21).

It suffices to show that each jump keeps the initial order. Consider the jump associated to a Poisson event at time $s \in N_{x,y}$ with uniform variable $U(s)$. There are two possibilities:
(1) If $x \in \Lambda(\lambda, \alpha)$, then $\lambda_x = \lambda_x^\alpha$. Since the function $g(\cdot)$ is monotone and the random variable $U(s)$ is the same for both marginals, the order is kept.
(2) If $x \in \Lambda^c(\lambda, \alpha)$, then $\lambda_x < \lambda_x^\alpha$. In this case a $\eta_{s-}(x)$ particle jumps from $x$ to $y$ if $U(s) < \lambda_x g(\eta_{s-}(x))$ and a $\eta_{s-}^\alpha(x)$ particle jumps from $x$ to $y$ if $U(s) < \lambda_x^\alpha g(\eta_{s-}^\alpha(x))$. Hence, if $\eta_{s-}(x) \leq \eta_{s-}^\alpha(x)$ and $\eta_{s-}(y) \leq \eta_{s-}^\alpha(y)$, then $\eta_s(y) \leq \eta_s^\alpha(y)$. On the other hand, $\eta_s(x) \leq \eta_s^\alpha(x) = \infty$. $\square$

To prove Proposition 1, we need the following result. It helps to prove that the second class particles do not stop forever at some place: eventually every such particle either move or coalesce.

**Lemma 2.** *Fix an environment $\lambda$ and consider the stationary process $(\eta_t^\alpha, t \in \mathbb{R})$ with time-marginal distribution $\nu_\lambda^\alpha$ and fix $x \in \Lambda(\lambda, \alpha)$. Then $\eta_t^\alpha(x) = 0$ infinitely often with probability one:*

$$\text{(22)} \qquad \liminf_{t \to \infty} \eta_t^\alpha(x) = 0.$$

*Proof.* Consider the discrete time stationary process $(\eta_n^\alpha(x), n \in \mathbb{N})$ —this is just the process $(\eta_t(x), t \in \mathbb{R})$ observed at integer times. It is sufficient to show

$$\text{(23)} \qquad \liminf_{n \to \infty} \eta_n^\alpha(x) = 0$$

with probability one. A theorem of Poincaré (Chapter IV in [21] or Theorem 3.4 of Chapter 6 in [6]) implies that for every $k \in \mathbb{N}$,

$$\mathbb{P}\big(\eta_n^\alpha(x) = k \text{ infinitely often in } n \mid \eta_0^\alpha(x) = k\big) = 1.$$

Returning for a moment to the continuous time process $\eta_t^\alpha$, if at time $t$ site $x$ has at least one particle, then one of the particles at $x$ will jump with probability bounded below by $g(1)\lambda_x/(1 + \lambda_x) > 0$, this is the probability the exponential jump time of $x$ is smaller than the jump-times of particles from the other sites to $x$, whose rate is bounded by $g(\infty) \sum_y p(y, x) = 1$. Fix $k \in \mathbb{N}$. By the same reasoning, for any $m$, if $\eta_m^\alpha(x) = k$, then there is a positive probability to be visiting 0 at time $m + 1$ independently of previous visits and uniformly in the configuration outside $x$ at time $m$. Since these are independent attempts, Borel–Cantelli implies

$$\mathbb{P}\big(\eta_n^\alpha(x) = 0 \text{ infinitely often in } n \mid \eta_0^\alpha(x) = k\big) = 1.$$

This implies (23). $\square$

*Proof of Proposition 1.* In an environment $\lambda$, consider the coupling process of two versions of the process $\eta_t^\alpha$ obtained by using the same family of Poisson processes $(N_{x,y} : x, y \in \mathbb{Z}^d)$ and uniform random variables $(U(s), s \in \bigcup_{x,y} N_{x,y})$. By $\{\bar{S}_\lambda^\alpha(t) : t \geq 0\}$ denote the semigroup of the process and by $\mathbb{P}$ the probability associated to the process.

Since $\nu_\lambda^\alpha$ is invariant for $\eta_t^\alpha$, it is enough to show that, for any $\alpha > 0$, any $\mu$ satisfying (3), any $\lambda$ ($\mathcal{P}$-a.s.), and for any $x \in \Lambda(\lambda, \alpha)$,

$$(24) \qquad \lim_{t \to \infty} (\mu^\alpha \times \nu_\lambda^\alpha) \bar{S}_\lambda^\alpha(t) \{(\xi, \eta) : \xi(x) \neq \eta(x)\} = 0.$$

In coupling terms, (24) reads

$$(25) \qquad \lim_{t \to \infty} \int \int \mu^\alpha(d\xi) \, \nu_\lambda^\alpha(d\eta) \, \mathbb{P}\big(\xi_t(x) \neq \eta_t(x) \mid (\xi_0, \eta_0) = (\xi, \eta)\big) = 0,$$

where we have denoted $\xi_t$ the first coordinate of the coupled processes and $\eta_t$ the second. Therefore, to prove the proposition it is enough to prove that, for any $\alpha > 0$, any $\mu$ satisfying (3), any $\lambda$ ($\mathcal{P}$-a.s.), any $\xi^0$ ($\mu$-a.s.), any $\eta^0$ ($\nu_\lambda^\alpha$-a.s.), and for any $x \in \Lambda(\lambda, \alpha)$,

$$(26) \qquad \lim_{t \to \infty} \mathbb{P}\big(\xi_t(x) \neq \eta_t(x) \mid (\xi_0, \eta_0) = (\xi^0, \eta^0)\big) = 0.$$

Without loss of generality we assume $x = 0$ and $\alpha$ small enough such that $0 \in \Lambda(\lambda, \alpha)$. Fix $\alpha$, $\lambda$, $\xi^0$ and $\eta^0$. The configurations $\xi^0$ and $\eta^0$ are in principle not ordered: there are sites $y \in \Lambda(\lambda, \alpha)$ such that $(\xi^0(y) - \eta^0(y))^+ > 0$ and sites $z \in \Lambda(\lambda, \alpha)$ such that $(\xi^0(z) - \eta^0(z))^- > 0$. We say that we have $\xi\eta$-discrepancies in the first case and $\eta\xi$-discrepancies in the second one.

Denote $\bar{\xi}_t(z) := \min\{\xi_t(z), \eta_t(z)\}$ the number of coupled particles at site $z$ at time $t$. The $\bar{\xi}$-particles move as regular zero-range particles; they are usually called *first class particles*. There is at most one type of discrepancy at each site at any time. Discrepancies of both types move as *second class particles*, i.e., $\xi\eta$-discrepancy jumps from $y$ to $z$ with rate

$$(27) \qquad \lambda_y^\alpha p(z - y)[g(\xi(y)) - g(\bar{\xi}(y))]$$

and $\eta\xi$-discrepancy jumps from $y$ to $z$ with rate

$$(28) \qquad \lambda_y^\alpha p(z - y)[g(\eta(y)) - g(\bar{\xi}(y))]$$

that is, second class particles jump with the difference rate. For instance, in the case $g(k) \equiv 1$, the second class particles jump only when there are no coupled particles in the site.

If a $\xi\eta$-discrepancy jumps to a site $z$ occupied by at least one $\eta\xi$-discrepancy, then the $\xi\eta$-discrepancy and one of the $\eta\xi$-discrepancies at $z$ coalesce into a coupled $\bar{\xi}$-particle in $z$. Analogously, for the case when a $\eta\xi$-discrepancy jumps to a site $z$ occupied by at least one $\xi\eta$-discrepancy. The coupled particle behaves from this moment on as a first class particle. If a discrepancy of any type jumps to a site $z$ with infinite number of particles, that is, $z \in \Lambda^c(\lambda, \alpha)$, then the discrepancy disappears. All particles in sites $x \in \Lambda^c(\lambda, \alpha)$ are first class $\bar{\xi}$-particles. Therefore, any particle that jump from any site $x \in \Lambda^c(\lambda, \alpha)$ is a first class particle.

At time zero there are $|\xi^0(y) - \eta^0(y)|$ discrepancies at site $y$. To the $i$th discrepancy at site $y$ at time zero, that is, discrepancy $(y, i)$, we associate the random walk $\zeta_n(y, i)$ from the model of Section 4.

Since the interaction with the other particles and the environment $\lambda$ governs the waiting times between jumps but does not affect the skeleton of the discrepancy motion until coalescence or absorbing time, it is possible to couple the skeleton of the discrepancy $(y, i)$ with the random walk $\zeta_n(y, i)$ in such a way that they perform the same jumps together until (a) the coalescence of the discrepancy with another discrepancy of different type or (b) the absorption of the discrepancy at some site of $\Lambda^c(\lambda)$. In any case, the number of discrete jumps is at most $\tau(y,i;\Lambda^c(\lambda))$. Therefore, the full trajectory of discrepancy $(y,i)$ is shorter (visits not more sites and has not more number of visits to each site) than the trajectory of the random walk $\zeta_n(y,i)$ in the time interval $[0, \tau(y,i;\Lambda^c(\lambda))]$. Thus, Proposition 2 implies that only a finite number of discrepancies visit $x$ and the number of visits of site $x$ by each of the discrepancies is finite.

Lemma 2 implies that there are no $\eta$-particles at $x$ infinitely often. Therefore, there are no $\eta\xi$-discrepancies at $x$ infinitely often. This means that every $\eta\xi$-discrepancy that at some moment is at $x$ will eventually jump out or coalesce. It follows that after some random time there is no $\eta\xi$-discrepancies in $x$ forever.

Moreover, if at time $t$ site $x \in \Lambda(\lambda, \alpha)$ has no $\eta$-particles, then a $\xi\eta$-discrepancy at $x$ will jump with probability bounded below by $g(1)\lambda_x/(1+\lambda_x) > 0$. Therefore, using Lemma 2, we see that after some random time there is no $\xi\eta$-discrepancies in $x$ forever. □

**Acknowledgments.** We thank Enrique Andjel and James Martin for discussions.

Part of this paper was written when the first author was participating of the program Principles of the Dynamics of Non-Equilibrium Systems in the Isaac Newton Institute for Mathematical Sciences visiting Newton Institute, in May-June 2006.

This paper was partially supported by FAPESP and CNPq.